\documentclass[a4paper]{article}
\usepackage[utf8]{inputenc}
\usepackage{amsmath,amssymb,amsthm}
\usepackage{graphicx}
\usepackage{url}
\usepackage{color}
\usepackage[realmainfile]{currfile}
\newtheorem{theorem}{Theorem}

\newtheorem{lemma}[theorem]{Lemma}

\numberwithin{theorem}{section}
\numberwithin{equation}{section}
\title{On the supremum of a quotient of power sums}

\author{
Stefan Gerhold \\
TU Wien \\
\tt{sgerhold@fam.tuwien.ac.at}
\and
Friedrich Hubalek \\
TU Wien \\
\tt{fhubalek@fam.tuwien.ac.at}
}

\date{\today}
\begin{document}

\maketitle

\begin{abstract}
  We define a function of two real vectors by a certain homogeneous quotient
  involving power sums, and show that its supremum grows
  asymptotically linearly w.r.t.\ the dimension.
  From this, we deduce a condition under which a parametric set of real matrices
  satisfies a set  of polynomial positivity constraints.
  This characterization finds an application in mathematical
  finance, in a recent study on price impact models.
\end{abstract}



\section{Introduction}

Reznick~\cite{Re83} has shown, among other results, that the maximum of 
\[
  \textstyle{\Big(\sum\limits_{i=1}^n x_i\Big)\Big(\sum\limits_{i=1}^n x_i^3\Big)\Big/\Big(\sum\limits_{i=1}^n x_i^2\Big)^2}
\]
over $\mathbb{R}^n$ is of order~$\sqrt{n}$, as $n\to\infty$. Our contribution is an analytical investigation of a homogeneous quotient involving {\em differences} of power sums in high dimension, motivated by recent work
on price impact models in mathematical finance~\cite{HeMaMu25}.
With $M_p(x):=\sum_{i=1}^n x_i^p$ for $p>0$,
 we prove the inequality
\begin{align*}
  Q(x,y) &:= \frac{\big(M_1(x)-M_1(y)\big)\big(M_2(y)-M_2(x)\big)}{M_3(x)+M_3(y)} \\
  &= \frac{\Big(\sum\limits_{i=1}^n x_i-\sum\limits_{i=1}^n y_i\Big)\Big(\sum\limits_{i=1}^n y_i^2-\sum\limits_{i=1}^n x_i^2\Big)}{\sum\limits_{i=1}^n x_i^3 + \sum\limits_{i=1}^n y_i^3}
  <   \frac{7\sqrt{7}-17}{27} n,
  \quad x,y \in \mathbb{R}^n_{>0},
\end{align*}
and show that the constant factor is sharp as $n\to\infty$.
Put differently, we maximize the product
$\big(\|x\|_1-\|y\|_1\big)\big(\|y\|_2^2-\|x\|_2^2\big)$
 over the $\ell^3$-unit ball, with a focus on large dimension.
 Note that this product is non-positive, if the two vectors $x,y\in\mathbb{R}^n_{>0}$ are ordered component-wise. It is easy to see, though, and will be shown below, that the supremum is always positive for $n\geq2$. 
 In what follows, we will also consider~$Q(x,y)$ for vectors
$x\in\mathbb{R}^{n+1}_{>0}$ and $y\in\mathbb{R}^n_{>0}$ of unequal length.
In Section~\ref{se:appl}, we discuss an application: 
The bound allows to characterize a certain parametric set of matrices that
arises in~\cite{HeMaMu25}.
This is also where the constraint
of considering only vectors with positive entries comes from.
As a side remark, we note that there is a considerable computational literature on optimization
of homogeneous polynomials, and refer to~\cite{HoSo14} and the references therein.

\section{Main result}

\begin{theorem}\label{thm:main}
  With
  \[
    c^* := \frac{7\sqrt{7}-17}{27} \approx 0.0563,
  \]
  we have
  \begin{equation}\label{eq:Q ineq}
     Q(x,y) < c^* n, \quad x,y \in\mathbb{R}^n_{>0},\ n\geq1,
  \end{equation}
  and
  \begin{equation}\label{eq:Q lim}
    \lim_{n\to\infty}  \frac1n \sup\big\{Q(x,y):x,y\in\mathbb{R}^n_{>0}\big\} = c^*.
  \end{equation}
  Both assertions remain true if the condition $x\in\mathbb{R}^n_{>0}$
  is replaced by $x\in\mathbb{R}^{n+1}_{>0}$.
\end{theorem}
First, we verify that~$Q$ has positive values for all $n\geq2$.
\begin{lemma}\label{le:pos}
  For $n=1$, we have $\max_{x,y>0} Q(x,y) = 0$. Furthermore,
  \[
     \sup\big\{Q(x,y):x,y\in\mathbb{R}^n_{>0}\big\} >0, \quad n\geq 2,
  \]
  and
  \[
     \sup\big\{Q(x,y):x\in\mathbb{R}^{n+1}_{>0},\ y\in\mathbb{R}^n_{>0}\big\} >0,
     \quad n\geq1.
  \]
\end{lemma}
\begin{proof}
  For $n=1$, note that
  \[
     Q(x,y) = \frac{(x-y)(y^2-x^2)}{x^3+y^3} = -\frac{(x-y)^2(x+y)}{x^3+y^3} \leq 0, \quad x,y>0.
  \]
  For $n\geq2$, define the vectors
  \[
    \hat{x}^{(n)} :=
    (\underbrace{1,\dots,1}_{\lceil n/2 \rceil},
      \underbrace{n^{-1},\dots,n^{-1}}_{\lfloor n/2 \rfloor}) \in \mathbb{R}^{n}_{>0},
  \]
  \[
    \hat{y}^{(n)} := \big(\tfrac{701}{1000}, \dots, \tfrac{701}{1000}\big)
    \in \mathbb{R}^{n}_{>0},
  \]
  and define $\tilde{x}^{(n)}\in\mathbb{R}^{n+1}_{>0}$ by adding another
  component $n^{-1}$ to $\hat{x}^{(n)}$.
  For even~$n$, $Q(\hat{x}^{(n)},\hat{y}^{(n)})$ is a rational function of~$n$.
  Using a computer algebra system, it is easily verified that this rational
  function is positive for real $n\geq 4$. In particular, $Q(\hat{x}^{(n)},\hat{y}^{(n)})>0$
  for integral $n\geq 4$. Analogously, the same holds for odd $n\geq 5$, and we have
  $Q(\tilde{x}^{(n)},\hat{y}^{(n)})>0$ for $n\geq7$.
  We can also use computer algebra
  to find  vectors $x,y$ with $Q(x,y)>0$ for the remaining finitely many cases.
  For instance, we have $Q(\check{x},\check{y})\approx 0.031$ for
  \begin{align*}
    \check{x} &:= \big( \tfrac32, \tfrac{1}{16\,777\,216}, \tfrac{1}{8\,388\,608},
    \tfrac{1}{2\,097\,152}, \tfrac{1}{8192}, \tfrac{1}{4096}, \tfrac34 \big)
    \in\mathbb{R}^7_{>0}, \\
    \check{y} &:= \big( 1, \tfrac{1}{256}, \tfrac{1}{64}, \tfrac{1}{16},
    \tfrac12, 1\big)  \in\mathbb{R}^6_{>0}. \qedhere
  \end{align*}
\end{proof}
\begin{proof}[Proof of Theorem~\ref{thm:main}]
  We begin with the upper estimate. For $n=1$, we have
  \[
    \sup_{x,y>0} Q(x,y) = 0< c^*
  \] 
   by Lemma~\ref{le:pos}.
  To prove the upper bound~\eqref{eq:Q ineq} for $n\geq 2$,  it is convenient to work on $\mathbb{R}^n_{\geq0}\setminus\{0\}$ instead of $\mathbb{R}_{>0}^n$. Suppose that~$Q$
  has a local maximum at $(x,y)=(\xi,\eta) \in \mathbb{R}_{\geq0}^{2n}\setminus\{(0,0)\}$.
  The global maximal value
  is positive by Lemma~\ref{le:pos},
  and so we may consider the logarithmic derivatives
  \[
    \partial_{x_i} \log Q(x,y) \big|_{(x,y)=(\xi,\eta)}  = U(\xi_i), \quad
    \partial_{y_i} \log Q(x,y) \big|_{(x,y)=(\xi,\eta)} = V(\eta_i),
  \]
  where
  \[
    U(t) := \frac{1}{M_1(\xi)-M_1(\eta)} - \frac{2t}{M_2(\eta)-M_2(\xi)}
      - \frac{3t^2}{M_3(\xi)+M_3(\eta)}
  \]
  and
  \[
    V(t) := {-\frac{1}{M_1(\xi)-M_1(\eta)}} +\frac{2t}{M_2(\eta)-M_2(\xi)}
      - \frac{3t^2}{M_3(\xi)+M_3(\eta)}.
  \]
  If $\xi_i>0$, then we must have $U(\xi_i)=0$, and if $\eta_i>0$, then $V(\eta_i)=0$.
  As each of the two quadratic polynomials $U,V$ has at most two zeros, we conclude that, for a
  global   maximum at $(\xi,\eta)$, their exist four numbers $\alpha,\beta,\gamma,\delta$ with
  \[
    \xi_i \in \{\alpha,\beta,0\}, \quad \eta_i \in \{\gamma,\delta,0\},
    \quad 1\leq i\leq n.
  \] 
  It thus suffices to consider
  \begin{equation}\label{eq:x y}
    x = (\underbrace{\alpha,\dots,\alpha}_{i}, \underbrace{\beta,\dots,\beta}_{j},0,\dots,0),
    \quad
    y = (\underbrace{\gamma,\dots,\gamma}_{k}, \underbrace{\delta,\dots,\delta}_{l},0,\dots,0)
  \end{equation}
  with
  \begin{align}
    \alpha > \beta > 0, & \quad \gamma > \delta >0, \notag \\
    1\leq i \leq n, &\quad 0\leq j \leq n-i, \label{eq:constraints} \\
    1\leq k \leq n, &\quad 0\leq l \leq n-k. \notag
  \end{align}
  Then, for $x,y$ as in~\eqref{eq:x y}, we have
  \begin{equation}\label{eq:Q S}
    Q(x,y) = \frac{S_1 S_2}{S_3},
  \end{equation}
  where
  \begin{align*}
    S_1 &:= i \alpha + j \beta -k\gamma - l \delta, \\
    S_2 &:= -i \alpha^2 - j \beta^2 +k\gamma^2 + l \delta^2, \\
    S_3 &:= i \alpha^3 + j \beta^3 +k\gamma^3 + l \delta^3.
  \end{align*}
  Thus, we arrive at a rational function of $(\alpha,\beta,\dots,k,l)$.
  By Lemma~\ref{le:pos}, the supremum of this function in the domain
  defined by~\eqref{eq:constraints}, with real variables $i,j,k,l$, is positive.
  By homogeneity and symmetry, we may  assume that~$\alpha$ and~$\gamma$
  are bounded.
  We proceed to determine the maximum of 
  \[
   h(\alpha,\beta,\dots,k,l):=\log (S_1 S_2/S_3),
  \]
  which yields an upper bound for $Q(x,y)$. 
  We define $h:=-\infty$ for $S_1 S_2\leq 0$, but note that~$h$ is finitely-valued 
  and smooth in a neighborhood of any global maximum.
  The partial derivatives of~$h$ are
  \[
    h_\alpha = u(\alpha)i,\quad h_\beta = u(\beta)j,\quad h_\gamma= v(\gamma)k, \quad
    h_\delta = v(\delta)l,
  \]
  where
  \[
    u(t) := \frac{1}{S_1}-\frac{2t}{S_2} - \frac{3t^2}{S_3}, \quad
    v(t) := {-\frac{1}{S_1}}+\frac{2t}{S_2} - \frac{3t^2}{S_3}.
  \]
  Moreover,
   \[
    h_i = U(\alpha),\quad h_j = U(\beta),\quad h_k= V(\gamma), \quad
    h_l = V(\delta),
  \]
  where
  \[
    U(t) = \frac{t}{S_1}-\frac{t^2}{S_2} - \frac{t^3}{S_3}, \quad
    V(t) = {-\frac{t}{S_1}}+\frac{t^2}{S_2} - \frac{t^3}{S_3}.
  \]
  We claim that $j=l=0$ holds at any global maximum of~$h$.
  The following argument rests on the identities
  \begin{equation}\label{eq:1st id}
    \frac{h_\alpha}{\alpha i} - \frac{h_i}{\alpha^2} - \frac{h_\beta}{\beta j}
      +\frac{h_j}{\beta^2} = {-\frac{2(\alpha-\beta)}{S_3}}
  \end{equation}
  and
  \begin{equation}\label{eq:2nd id}
    \frac{h_\gamma}{\gamma k} - \frac{h_k}{\gamma^2} - \frac{h_\delta}{\delta l}
      +\frac{h_l}{\delta^2} = {-\frac{2(\gamma-\delta)}{S_3}}.
  \end{equation}
  They show that no point in the interior of the domain of~$h$ at
  which~$h$ is $>-\infty$ satisfies the
  first order conditions $h_\alpha = \dots = h_l=0$.
  Furthermore, if there was a maximum at the boundary $j=n-i$, then~\eqref{eq:1st id}
  would imply
  \[
    \frac{h_j}{\beta^2} = {-\frac{2(\alpha-\beta)}{S_3}}<0,
  \]
  but then~$h$ would \emph{decrease} as we cross the boundary in the $j$-direction.
  The only remaining possibility is a maximum at the boundary~$j=0$ (equivalently, we could put $\beta=0$),
  and similarly we can show from~\eqref{eq:2nd id} that $l=0$.
  Moreover, by symmetry, it suffices to consider $i\leq k$.
  For $n=2$, \eqref{eq:Q ineq} is now easily verified.
  By induction, we may assume
  $k=n$, because otherwise the parameters $\alpha,\gamma,i,k$ would
  be admissible for the problem with~$n$ replaced by~$k$, and we can simply
  note that $c^*k<c^*n$.
  Assuming $1=\alpha>\gamma$, by symmetry and homogeneity, we have formally reduced the problem of optimizing~$Q$ to choosing
  \[
    x= (\underbrace{1,\dots,1}_{i},0,\dots,0), \quad
    y = (\gamma,\dots,\gamma)
  \]
  with $1\leq i\leq n$ and $0<\gamma<1$.
  We will see, though, that the optimal~$i$ for the function~$h$
  is not an integer, which shows that the bound~\eqref{eq:Q ineq} is strict.
  We put $i=pn$ with unknown $p\in(0,1)$, and so our target becomes
  \[
    \frac{S_1 S_2}{S_3} = \frac{(i-n\gamma)(-i+n\gamma^2)}{i+n \gamma^3}
    =\frac{(p-\gamma)(-p+\gamma^2)}{p+\gamma^3} n.
  \]
  Maximizing the latter fraction is straightforward, as the system $\partial_p(..)=0$,
  $\partial_\gamma(..)=0$ can be solved explicitly, with a single positive solution
  \[
      p^* := \frac{16-5\sqrt7}{27} \approx 0.102,\quad \gamma^* := \frac{\sqrt7-2}{3}
      \approx 0.215.
  \]
  It is elementary to show that $(p-\gamma)(-p+\gamma^2)/(p+\gamma^3)$, 
  $(p,\gamma)\in(0,1)^2$, has a global
  maximum at $(p^*,\gamma^*)$, with value $c^*$.
  This finishes the proof of~\eqref{eq:Q ineq}.
  
  Inspired by the upper estimate, we use the vectors
  \begin{equation}\label{eq:x^n}
    x^{(n)} := (\underbrace{1,\dots,1}_{\lfloor p^* n\rfloor},n^{-1},\dots,n^{-1}) \in \mathbb{R}^n_{>0}
  \end{equation}
  and
  \begin{equation}\label{eq:y^n}
    y^{(n)} := (\gamma^*,\dots,\gamma^*) \in \mathbb{R}^n_{>0}
  \end{equation}
  to prove the required asymptotic lower bound. From
  \[
    M_\nu\big(x^{(n)}\big) = p^* n + \mathrm{O}(1), \quad \nu=1,2,3,
  \]
  as $n\to\infty$ and
  \[
    M_\nu\big(y^{(n)}\big) = (\gamma^*)^\nu n + \mathrm{O}(1), \quad \nu=1,2,3,
  \]
  we calculate
  \begin{align*}
     Q\big(x^{(n)},y^{(n)}\big) &=
     \frac{\big(p^*n - \gamma^* n + \mathrm{O}(1)\big)
       \big({-p^*n}+(\gamma^*)^2n  +\mathrm{O}(1)\big)}{p^* n + (\gamma^*)^3n + \mathrm{O}(1)} \\
     &\sim \frac{(p^*-\gamma^*)\big({-p^*}+(\gamma^*)^2\big)}{p^*+(\gamma^*)^3}n = c^* n,
     \quad n\to\infty.
  \end{align*}
  The proof of~\eqref{eq:Q ineq} and~\eqref{eq:Q lim} is finished. Proving the remaining
  assertions is a simple modification: The vector~\eqref{eq:x^n}
  from the lower bound receives an additional component~$n^{-1}$.
  The range of~$i$ in the upper bound becomes $1\leq i\leq n+1$, which makes no essential difference.
  For the induction base, we need to compute
  \[
    \sup\big\{ Q(x,y) : x\in\mathbb{R}^2_{>0},\ y \in\mathbb{R}_{>0} \big\}.
  \]
  We may assume~$x_1=1$, and it is then an easy exercise that the maximum of
  \[
    Q(x,y) = \frac{(1+x_2-y)(y^2-1-x_2^2)}{1+x_2^3+y^3}
  \]
  is located at $x_2=1$, $y=\tfrac12(\sqrt{9+4\sqrt{6}}-1)$, and has a value
  that is smaller than~$2c^*$.
\end{proof}

\section{Application}\label{se:appl}

For $d\geq 2$, we define the set~$\mathcal{M}_d$ of matrices $M=(m_{ij}) \in \mathbb{R}^{d\times d}$
such that
\begin{equation}\label{eq:Psi}
  \Psi_M(z,s) := \sum_{l=1}^d\Big( m_{ll} z_l^3 + s_l z_l \sum_{k\neq l}m_{lk}s_k z_k^2 \Big) >0
\end{equation}
for all $s\in\{\pm1\}^d$ and $z\in\mathbb{R}^d_{>0}$.
In general, this gives $2^{d-1}-1$ inequalities which have to be
satisfied by all  $z\in\mathbb{R}^d_{>0}$, as~$s$ and~$-s$
yield the same condition, and $s=(1,\dots,1)$ can of course be discarded.
\begin{lemma}\label{le:cone}
  The set~$\mathcal{M}_d$ is a convex cone. Any matrix
  in~$\mathcal{M}_d$ has non-negative diagonal elements.
  If $M$ is diagonally dominant in the sense that
  \begin{equation}\label{eq:diag dom}
    m_{ii}> \sum_{l=1}^d \sum_{k\neq l}|m_{lk}|, \quad 1\leq i\leq d,
  \end{equation}
  i.e.\ each diagonal element dominates the sum of the absolute values of all
  off-diagonal elements, then $M\in \mathcal M_d$. In particular,
  diagonal matrices with positive diagonal elements are in $\mathcal M_d$.
\end{lemma}
\begin{proof}
  The first statement is obvious. The second one is also
  clear, by considering large multiples of unit vectors. Now let $z\in \mathbb{R}^d_{>0}$ and suppose that $z_i$ is one of the maximal entries of the vector. Then we have  for all $s\in\{\pm1\}^d$
  \begin{align*}
    \Psi_M(z,s) &= \sum_{l=1}^d\Big( m_{ll} z_l^3 + s_l z_l \sum_{k\neq l}m_{lk}s_k z_k^2 \Big) \\
    &\geq m_{ii} z_i^3 - z_i^3  \sum_{l=1}^d\sum_{k\neq l}|m_{lk}| >0. \qedhere
  \end{align*}
\end{proof}
The set~$\mathcal{M}_d$ arises naturally in~\cite{HeMaMu25}. As elaborated there, positivity of the term 
$h(\boldsymbol{J}_t)^\top \boldsymbol{\zeta} \boldsymbol{J}_t$ in~(6.1)
of~\cite{HeMaMu25}, for any vector $\boldsymbol{J}_t=(J_t^i)$, is crucial for the underlying model to make sense. For the standard
choice $h(x)=\mathrm{sgn}(x)\sqrt{|x|}$ (square-root price impact), and with
$z_i=\sqrt{|J^i_t|}$ and $M=\boldsymbol{\zeta}$, this amounts to the condition~\eqref{eq:Psi}.
It is thus of interest to find conditions
that ensure $M\in\mathcal{M}_d$, in particular for large dimensions.
Presumably, it is hard to give a simple general if-and-only-if condition.
Even for $d=2$, where a full characterization is possible, a complicated case distinction
involving algebraic expressions of the matrix elements $m_{ij}$ arises.
Beyond the trivial case of diagonal matrices of arbitrary dimension, a reasonable next step is to consider equal off-diagonal elements,
and so we define
\begin{equation}\label{eq:def Mb}
    M_d(b) := \begin{pmatrix}
     1 & b & \dots & b \\
     b & 1 & & \vdots \\
     \vdots & & \ddots & b \\
     b & \cdots & b & 1
     
    \end{pmatrix},
    \quad b\in \mathbb R_{> 0},
\end{equation}
for which~\eqref{eq:Psi} reduces to the condition
\begin{equation}\label{eq:Psi spec}
  \sum_{l=1}^dz_l^3 + b \sum_{l=1}^d s_l z_l \sum_{k\neq l}s_k z_k^2  >0,
  \quad s\in\{\pm1\}^d,\ z\in\mathbb{R}^d_{>0}.
\end{equation}
\begin{lemma}\label{le:b}
  For $d\geq2$ and $b>0$, \eqref{eq:Psi spec} is equivalent to
  \[
    \frac1b > 1+ Q(x,y), \quad x\in\mathbb{R}^{\lceil d/2 \rceil}_{>0},\
    y \in \mathbb{R}^{\lfloor d/2 \rfloor}_{>0}.
  \]
\end{lemma}
\begin{proof}
  We assume that~$d$ is even, as the proof for odd~$d$
  is analogous.
  By symmetry of~\eqref{eq:Psi spec}, only the number of positive entries in~$s\in\{\pm1\}^d$
  is relevant, not their position. Moreover, it suffices to consider
  \begin{equation}\label{eq:s}
    s=(\underbrace{-1,\dots,-1}_{d}, \underbrace{1,\dots,1}_{d}).
  \end{equation}
  Indeed, the other sign combinations~$s$ lead to a larger
   number of $(+1)$-s in the matrix $(s_i s_j)$. As the $z_i$ are positive, this 
  makes the desired inequality~\eqref{eq:Psi spec} easier to satisfy.
  For~$s$ as in~\eqref{eq:s},
  the second sum in~\eqref{eq:Psi spec} becomes
  \begin{align*}
    \sum_{l=1}^d s_l z_l & \sum_{k\neq l}s_k z_k^2
      = \sum_{l=1}^{d/2} s_l z_l \Bigg(\sum_{\substack{k=1 \\ k\neq l}}^{d/2} s_k z_k^2
       + \sum_{k=d/2+1}^d s_k z_k^2\Bigg) \\
    &\quad \qquad \qquad \qquad + \sum_{l=d/2+1}^{d} s_l z_l \Bigg(\sum_{k=1 }^{d/2} s_k z_k^2
       + \sum_{\substack{k=d/2+1\\k\neq l}}^d s_k z_k^2\Bigg) \\
   &= \sum_{l=1}^{d/2}z_l \sum_{\substack{k=1 \\ k\neq l}}^{d/2} z_k^2
      -  \sum_{l=1}^{d/2}z_l \sum_{k=d/2+1}^d z_k^2
      - \sum_{l=d/2+1}^d z_l \sum_{k=1}^{d/2} z_k^2
      + \sum_{l=d/2+1}^{d}z_l \sum_{\substack{k=d/2+1 \\ k\neq l}}^{d} z_k^2 \\
&=  \sum_{l=1}^{d/2}z_l\bigg( \sum_{k=1}^{d/2}z_k^2 - z_l^2\bigg)
    - \sum_{l=1}^{d/2}z_l  \sum_{k=d/2+1}^d z_k^2 \\
    & \qquad \qquad \qquad \quad
    - \sum_{l=d/2+1}^d z_l \sum_{k=1}^{d/2}z_k^2
    +\sum_{l=d/2+1}^d z_l \bigg(\sum_{k=d/2+1}^d z_k^2 -z_l^2 \bigg) \\
&= -M_3(x)- M_3(y)-\big( M_1(x)-M_1(y)\big) \big(M_2(y)-M_2(x)\big),
  \end{align*}
  where
    \[
    x := (z_1,\dots,z_{d/2}), \quad
    y := (z_{d/2+1},\dots,z_{d}).
  \]
  As the first sum in~\eqref{eq:Psi spec} equals $M_3(x)+M_3(y)$, the statement follows.
\end{proof}
We can now state our main result on membership of~\eqref{eq:def Mb}
in the cone~$\mathcal{M}_d$. While the naive criterion~\eqref{eq:diag dom}
shows that $b\lesssim 1/d^2$ suffices, we can infer from Theorem~\ref{thm:main}
that the true asymptotic order of the
maximal value of~$b$ is $1/d$ for large~$d$.

\begin{theorem}\label{thm:main b}
  For $d\geq 2$, the set of admissible~$b$ for the matrices~\eqref{eq:def Mb} is an interval with closure
  \begin{equation}\label{eq:M bd}
    \mathrm{cl}\{ b > 0 : M_d(b) \in \mathcal{M}_d \} = [0,b_d],
  \end{equation}
  where
  \begin{equation}\label{eq:bd sup}
    b_d := \frac{1}{1+\sup\Big\{Q(x,y): x\in\mathbb{R}^{\lceil d/2 \rceil}_{>0},\
     y \in \mathbb{R}^{\lfloor d/2 \rfloor}_{>0} \Big\}}.
  \end{equation}
  The sequence $b_d>0$ decreases, satisfies   the lower bound
   \begin{equation}\label{eq:bd lower}
    b_d \geq \frac{1}{1+c^* \lfloor d/2 \rfloor}, \quad d \in \mathbb N,
  \end{equation}
  and has the asymptotics
   \begin{equation}\label{eq:bd asympt}
    b_d \sim \frac{2}{c^*d}=\frac{54}{7\sqrt{7}-17}\ \frac1d
    \approx \frac{35.52}{d},\quad d\to\infty.
   \end{equation}
\end{theorem}
\begin{proof}
  It follows from Lemma~\ref{le:cone} that the set $\{ b > 0 : M_d(b) \in \mathcal{M}_d \}$
  is a non-empty interval.
  By Lemmas~\ref{le:pos} and~\ref{le:b}, we have  $M_d(b) \in \mathcal{M}_d$ if and only if
  \[
    b < \frac{1}{1+ \max\{Q(x,y),0\}}, \quad x\in\mathbb{R}^{\lceil d/2 \rceil}_{>0},\
    y \in \mathbb{R}^{\lfloor d/2 \rfloor}_{>0}.
  \]
  This proves~\eqref{eq:M bd}, and~\eqref{eq:bd lower} and~\eqref{eq:bd asympt}
  follow from Theorem~\ref{thm:main}.
  To see that~$b_d$ decreases, let $d\geq 3$ and $b,\varepsilon>0$ be such that
  $M_d(b(1+\varepsilon)) \in \mathcal{M}_{d}$. Thus, Lemma~\ref{le:b} implies
  \[
    1 > b(1+\varepsilon)\big(1+Q(x,y)\big),\quad x\in\mathbb{R}^{\lceil d/2 \rceil}_{>0},\
    y \in \mathbb{R}^{\lfloor d/2 \rfloor}_{>0}.
  \]
  Taking $x_{\lceil d/2 \rceil} \downarrow 0$ for odd~$d$ respectively
  $y_{ d/2 } \downarrow 0$ for even~$d$  yields
  \[
    1 \geq b(1+\varepsilon)\big(1+Q(x,y)\big)
    > b \big(1+Q(x,y)\big),
    \quad x\in\mathbb{R}^{\lceil (d-1)/2 \rceil}_{>0},\
    y \in \mathbb{R}^{\lfloor (d-1)/2 \rfloor}_{>0}.
  \]
  Therefore, $M_{d-1}(b)\in\mathcal{M}_{d-1}$. As $\varepsilon>0$ was arbitrary,
  we have shown that $b_{d-1}\geq b_d$.
\end{proof}

\section{Numerical values}

For $d=3$, the optimization problem in~\eqref{eq:bd sup} can be solved
explicitly, as stated at the end of the proof of Theorem~\ref{thm:main}. The number $b_3$ is the only positive root of the polynomial
\begin{equation}\label{eq:4poly}
  20x^4 + 60x^3 + 9 x^2 -54 x -27,
\end{equation}
with explicit expression
\begin{equation}\label{eq:4expl}
  b_3 = 
  \tfrac14\Big(\sqrt{\tfrac35 (39+16 \sqrt{6})}-3\Big)
  \approx 0.962.
\end{equation}
Alternatively, $b_d$ can be computed by quantifier elimination~\cite{St17} 
from~\eqref{eq:Psi spec}, which also yields~\eqref{eq:4poly} and~\eqref{eq:4expl}. As mentioned in the proof of Lemma~\ref{le:b},
it suffices to consider~$s$ as in~\eqref{eq:s} for even~$d$.
Analogously, for odd~$d$, the inequality~\eqref{eq:Psi spec} needs to be checked only for
\[
  s=(\underbrace{-1,\dots,-1}_{(d+1)/2},\underbrace{1,\dots,1}_{(d-1)/2}).
\]
Still, quantifier elimination is computationally hard, and so we could calculate~$b_d$
with this method only until $d=4$. It turns out that $b_3=b_4$. In principle, upper bounds can be computed
by using numerical optimization to approximate the supremum in~\eqref{eq:bd sup}.
It is not trivial, though, to find good starting values.
Note that the vectors~\eqref{eq:x^n} and~\eqref{eq:y^n} are not useful in that
respect for small~$n$. In fact, $Q(x^{(n)},y^{(n)})$ is negative for $n\leq9$.
To compute the values~$b_5$ and~$b_6$ in Table~1, we instead used Mathematica's {\tt FindInstance}
and {\tt Reduce} commands to show that~$Q$ in~\eqref{eq:bd sup} can be larger than $\frac{1079}{10000}$,
but not $\frac{1080}{10000}$.

Table~2 gives estimates for~$b_d$ for large~$d$.
For the upper estimates, we evaluate~$Q$ at the vectors~\eqref{eq:x^n}
with $n=\lceil d/2 \rceil$ and~\eqref{eq:y^n}
with $n=\lfloor d/2 \rfloor$, which yields an upper bound
for the right hand side of~\eqref{eq:bd sup}.
For large~$d$ our computations show that~$b_d$ is very close to $b_{d+1}$
for even~$d$. This is in contrast to the fact that $b_3=b_4$, as mentioned above.

\begin{table}[h]
\centering
\begin{tabular}{ c | c c c c c }
 $d$ & 2 & 3 & 4 &  5 & 6  \\ 
 \hline
 \text{Lower estimate} & 0.946 & 0.946 & 0.898 &  0.898 & 0.855 \\  
  $b_d$ & 1 & 0.962 &  0.962 & 0.902 & 0.902 
\end{tabular}
\caption{The numbers~$b_d$ from~\eqref{eq:bd sup} for small~$d$, and the
lower estimate~\eqref{eq:bd lower}.}
\end{table}
\begin{table}[h]
\centering
\begin{tabular}{ c | c c c c c c c }
 $d$ & 50 & 100 & 150 & 200 & 300 & 400 & 500 \\ 
 \hline
 \text{Lower estimate} & 0.415 & 0.262 & 0.191  & 0.150   & 0.105 & 0.081 & 0.066 \\  
 \text{Upper estimate} & 0.510 & 0.295  & 0.210  & 0.161  & 0.111 & 0.084 & 0.068 \\
 \text{Asymptotics} &  0.710 & 0.355  & 0.236 & 0.177  & 0.118 & 0.088 & 0.071
\end{tabular}
\caption{Estimates for~$b_d$ for larger~$d$. Note that~\eqref{eq:bd lower}, which
is asymptotically equal to~\eqref{eq:bd asympt}, gives a better
approximation than~\eqref{eq:bd asympt}.
 }
\end{table}

\bigskip

{\bf Acknowledgement:} We thank Johannes Muhle-Karbe for suggesting the problem
of studying the cone~$\mathcal{M}_d$.

\bibliography{literature}
\end{document}